\documentclass[11pt,a4paper,reqno]{amsart}
\usepackage{amsmath}
\usepackage{amssymb,latexsym}
\usepackage[dvips]{graphicx}
\usepackage{color}
\usepackage{cite}
\usepackage{amsthm}

\newcommand{\R}{\mathbb R}

\newtheorem{theorem}{Theorem} [section]
\newtheorem{lemma}{Lemma} [section]
\newtheorem{proposition}{Proposition} [section]
\newtheorem{corollary}{Corollary} [section]

\newtheorem{remark}{Remark}[section]

\let\ssection=\section\renewcommand{\section}{\setcounter{equation}{0}\ssection}
\begin{document}
\date{}

\address{M. Darwich: Universit\'e Fran\c{c}ois rabelais de Tours, Laboratoire de Math\'ematiques
et Physique Th\'eorique, UMR-CNRS 6083, Parc de Grandmont, 37200
Tours, France} \email{Mohamad.Darwich@lmpt.univ-tours.fr}
\title[Blowup]{On the $L^{2}$-critical nonlinear Schr\"odinger Equation with a nonlinear damping.}
\author{Darwich Mohamad.}

\begin{abstract}
We consider the Cauchy problem for the  $L^{2}$-critical  nonlinear Schr\"{o}dinger equation with a nonlinear damping.
 According to the power of the damping term, we prove  the global existence  or the
existence of finite time blowup dynamics with the log-log blow-up speed  for $\|\nabla u(t)\|_{L^2}$.
\end{abstract}

\maketitle
\section{Introduction}
In this paper, we study  the blowup and the global existence of solutions for the focusing NLS equation with a nonlinear damping $(NLS_{ap})$:
\begin{equation}
\left\{
\begin{array}{l}
iu_{t} + \Delta{u} +|u|^{\frac{4}{d}}u + ia|u|^pu =0,  (t,x) \in [0,\infty[\times \mathbb{R}^{d}, d=1,2,3,4. \\
u(0)= u_{0} \in H^1(\mathbb{R}^{d})
\end{array}
\right. \label{NLSap}
\end{equation}
with initial data $u(0)= u_{0} \in H^1(\mathbb{R}^{d})$  where $a > 0$ is the coefficient of friction and $ p\ge 1 $.
Note that  if we replace  $+|u|^{\frac{4}{d}}u$ by $-|u|^{\frac{4}{d}}u$ , (\ref{NLSap}) becomes the 
 defocusing NLS equation.\\

Equation (\ref{NLSap}) arises in various areas of nonlinear optics, plasma physics and fluid mechanics.
Fibich \cite{Fibich}  noted  that in the nonlinear optics context, the origin of the nonlinear
damping is multiphoton absorption. For example, in the case of solids the number $p$ corresponds
to the number of photons it takes to make a transition from the valence band to
the conduction band. Similar behavior can occur with free atoms, in this case $p$
corresponds to the number of photons needed to make a transition from the ground
state to some excited state or to the continuum.

 The Cauchy problem for  (\ref{NLSap}) was studied by  Kato \cite{Kato} and  Cazenave\cite{Cazenave1} and it is known that if $p < \frac{4}{d-2}$,
 then the problem is locally well-posed in $H^1(\mathbb{R}^{d})$: 
For any $u_{0} \in H^{1}(\mathbb{R}^{d})$, there exist $T \in (0,\infty]$ and a unique solution $u(t)$ of $(1.1)$
 with $u(0)=u_{0}$ such that $u \in C([[0,T]);H^1(\mathbb{R}^{d}))$. Moreover, T is the maximal existence time of the solution $u(t)$ in 
the sense that if $T < \infty$ then $\displaystyle{ \lim_{t\rightarrow T}{\|u(t)\|_{H^1(\mathbb{R}^{d})}}}=\infty$.\\
Let us notice that for $a=0$ (\ref{NLSap}) becomes the $L^2$-critical nonlinear Schr\"{o}dinger equation:\\
\begin{equation}
\left\{
\begin{array}{l}
 iu_{t} + \Delta u + |u|^{\frac{4}{d}}u = 0\\
 u(0)=u_{0} \in H^{1}(\mathbb{R}^{d})
 \end{array}
\right.\label{NLS}
 \end{equation}
  For $u_0 \in  H^1$, a sharp criterion for global existence for (\ref{NLS}) has been exhibited by Weinstein
\cite{weinst}:  Let  $Q$  be the unique positive solution to 
 \begin{equation}\label{ellip}
 \Delta Q + Q|Q|^{\frac{4}{d}} = Q.
 \end{equation}
For $\|u_0\|_{L^2} < \|Q\|_{L^2}$ , the solution of (\ref{NLS}) is global in $H^1$. This follows from the conservation
of the energy and the $L^2$ norm and the sharp Gagliardo-Nirenberg inequality:
\begin{equation}
\forall u \in H^1,  E(u) \geq \frac{1}{2}(\int |\nabla u|^2)\bigg(1 - \big(\frac{\int |u|^2}{\int |Q|^2}\big)^{\frac{2}{d}}\bigg). \nonumber
\end{equation} 
On the other hand, there exists explicit solutions with $\|u_0\|_{L^2} =\|Q\|_{L^2}$ that blow up in finite time in the regime $\frac{1}{T-t}$.\\               
In the series of papers \cite{MerleRaphael1,Merle6}, Merle and Raphael have studied the blowup for (\ref{NLS}) with
$ \|Q\|_{L^2}<  \|u_0\|_{L^2} < \|Q\|_{L^2} + \delta$, $\delta$ small
 and have proven the
  existence of the blowup regime corresponding to the log-log law:
\begin{equation}\label{speed}
\displaystyle{\|u(t)\|_{H^1(\mathbb{R}^{d})} \sim \bigg(\frac{\text{log}\left|\text{log}(T-t)\right|}{T-t}\bigg)^{\frac{1}{2}}.}
\end{equation} 
In  \cite{Darwich}, Darwich  has proved in  case of the linear damping  ($p=0$), the global existence  in $H^1$ for
$\|u_0\|_{L^2} \leq \|Q\|_{L^2}$, and has showed that the log-log regime is stable by  such perturbations (i.e. there exist  solutions 
 blows up in finite time with the log-log law).\\
 Numerical observations suggest that this finite time blowup phenomena persists in the case of the nonlinear
 damping for $p < \frac{4}{d}$ ( see Fibich \cite{Fibich} and \cite{passota}). Passot and  Sulem \cite{passota} have proved that 
the solutions  are global in $H^1(\mathbb{R}^{2})$ in the case where the power of the damping term is strictly greater to the focusing nonlinearity.
 The case where the power of the damping term is equal to the focusing nonlinearity,  
 "small damping prevents blow-up ? "was an open question for  Sparber and Antonelli in their paper \cite{Sparber} and for 
Fibich and Klein in their paper \cite{Fibich1}. Our results can gives  an answer for their open problem, at least for  the $L^2$-critical case.
In fact, our aim in 
this paper  is study for  each value of $(d,p)$,  the existence of blow-up solutions as well as  global existence criteria. 
And know if the regime log-log still stable by such perturbations.\\
\vspace{0,3cm}
Let us now our results:

\begin{theorem}\label{theoremprincipal}
Let  $u_0$ in $H^1(\mathbb{R}^{d})$ with $d=1,2,3,4$:
\begin{enumerate}
\item if $\frac{4}{d-2} > p \geq \frac{4}{d}$, then the solution of (\ref{NLSap}) is global in $H^{1}$.\vspace{2mm}
\item if $ 1\le p < \frac{4}{d}$ and $ 1\le p\le 2 $,  then there exists $ 0<\alpha <\| Q\|_{L^2} $ such that  for any $u_0 \in H^{1}$  with $\|u_0\|_{L^2} < \alpha$,  the emanating solution is global in $H^{1}$.
\vspace{2mm}
\item if $ 1\le p < \frac{4}{d}$, then there exists $\delta_0 > 0$ such that $\forall a > 0$ and $\forall \delta \in]0, \delta_{0}[$,
 there exists $u_0 \in H^1$ with $\|u_0\|_{L^2} = \|Q\|_{L^2} + \delta$, such that the solution of (\ref{NLSap}) blows up in finite time in the log-log regime.
\item if $ 1 \le p < \frac{4}{d}$, then
 there does not exists an intial data $u_0$ with $\|u_0\|_{L^2} \leq \|Q\|_{L^2}$ such that the solution $u$ of (\ref{NLSap})
 blow up in finite time with this law: 
$$\frac{1}{(T-t)^{\beta - \epsilon }} \lesssim \|\nabla u(t)\|_{L^2(\mathbb{R}^{d})} \lesssim \frac{1}{(T-t)^{\beta + \epsilon}},
$$
$\text{for any}~ \beta\in ]0,\frac{2}{pd}[~~ \text{and} ~~0 < \epsilon < \frac{2-\beta pd}{8+pd}$.
\end{enumerate}
\end{theorem}
\begin{remark}
 Note that, part (4) of Theorem \ref{theoremprincipal} prove in particular that
 we dont have the blowup in the regime log-log for any $p \in [1,\frac{4}{d}[$ and in the regime $\frac{1}{t}$ for 
$d=1$ and $1\leq p < 2$, for initial data with 
 critical or subcritical mass.
\end{remark}

In the "critical" case $p=\frac{4}{d}$, we have more precisely  :
\begin{theorem}\label{theorem2}
 Let $p = \frac{4}{d}$, then the initial-value probem (\ref{NLSap}) is globally well posed in $H^{s}(R^d)$, $s \geq 0$.
Moreover, there exist unique $u_+$ in $ L^2$ such that 
\begin{equation}\label{scattering}
 \|u(.,t)-e^{it\Delta}u_{+}\|_{L^2} \longrightarrow 0, ~ t \longrightarrow + \infty,
\end{equation}
where $e^{it\Delta}$ is the free evolution.
\end{theorem}

\begin{remark}
 Theorem \ref{theorem2}  and part (1) and (2) of Theorem \ref{theoremprincipal} still hold in  the defocusing case.
\end{remark}
\begin{remark}
 Note that if $u(t,.)$ is a solution of $(NLS_{a,p})$ then $\overline{u(-t,.)}$ is a solution of $(NLS_{-a,p})$, then we dont have the scattering in $-\infty$, because
 this changes the sign of  the coefficient of friction.
\end{remark}

$\textbf{Acknowledgments.}$ I  would like  to thank my advisor Prof. Luc Molinet for his rigorous attention to this paper,  Dr. Christof Sparber 
for his remarks and Prof. Baoxiang Wang for having given me the reference of Lemma \ref{Friedman}.
\section{Proof of part (4) of Theorem \ref{theoremprincipal}}

Special solutions play a fundamental role for the description of the dynamics of (NLS). They are the solitary waves of the form $u(t, x) =\exp(it)Q(x)$, where $Q$ solves:
 \begin{equation}\label{ellip}
 \Delta Q + Q|Q|^{\frac{4}{d}} = Q.
 \end{equation}
 The
pseudo-conformal transformation applied to the stationary solution $e^{it}Q(x)$
yields an explicit solution for (NLS)
$$
S(t,x) = \frac{1}{\mid t \mid^{\frac{d}{2}}} Q(\frac{x}{t})e^{-i \frac{\mid x \mid^2}{4t} + \frac{i}{t}}
$$
which blows up at $T = 0$.\\
 Note that 
\begin{equation}\label{nablaS}
\|S(t)\|_{L^2} = \|Q\|_{L^2} ~~\text{and}~~ \|\nabla S(t)\|_{L^2}\sim \frac{1}{t}
\end{equation}
 It turns out that
$S(t)$ is the unique minimal mass blow-up solution in $H^1$  up to the symmetries of the equation ( see \cite{Merleseul}).

A known lower bound ( see \cite{Merle4}) on the blow-up rate for (NLS) is 
\begin{equation}\label{lowerNLSap}
\|\nabla u(t)\|_{L^2} \geq \frac{C(u_0)}{\sqrt{T-t}}.
\end{equation}

Note that this blow-up rate is the one of $S(t)$ given by (\ref{nablaS}) and log-log given by (\ref{speed}). Now, we will prove the part (4) of Theorem \ref{theoremprincipal}. 
\\ 

 For this we need the following Theorem ( see \cite{Hmidi}) :\\
\begin{theorem}\label{limsup}
 Let $(v_{n})_{n}$ be a bounded family of $H^1(\mathbb{R}^{d})$, such that:
 \begin{equation}
 \limsup_{n \rightarrow +\infty}\left\|\nabla v_{n}\right\|_{L^2(\mathbb{R}^{d})} \leq M \quad and \quad \limsup_{n \rightarrow +\infty}\left\|v_{n}\right\|_{L^{\frac{4}{d} + 2}} \geq m.
 \end{equation}
 Then, there exists $(x_{n})_{n} \subset \mathbb{R}^{d}$ such that:
 \begin{equation}
 v_{n}(\cdot + x_{n}) \rightharpoonup V \quad weakly, \nonumber
 \end{equation}
 with $\left\|V\right\|_{L^2(\mathbb{R}^{d})} \geq (\frac{d}{d+4})^{\frac{d}{4}}\frac{m^{\frac{d}{2}+1} + 1}{M^{\frac{d}{2}}}\left\|Q\right\|_{L^2(\mathbb{R}^{d})}$.
 \end{theorem}
Let us recall the following quantities:\\
\\
$L^2$-norm : $\left\|u(t, x)\right\|_{L^2} = \displaystyle{\int |u(t,x) |^2dx}$.\\
Energy : $E(u(t, x)) = \frac{1}{2}\|\nabla u\|_{L^2}^{2} - \frac{d}{4 + 2d}\|u\|_{L^{\frac{4}{d}+2}}^{\frac{4}{d}+2}. $\\
Kinetic momentum : $P(u(t))=Im(\displaystyle{\int} \nabla u \overline{u}(t,x)).$\\
\begin{remark}
It is easy to prove that if $u$ is a solution of (\ref{NLSap}) on $[0,T[ $, then for all $ t\in [0,T[ $ it holds
 \begin{equation}\label{mass a}
 \frac{d}{dt}\|u(t)\|_{L^2}= -2a\int|u|^{p+2} \, , 
 \end{equation}
 \begin{equation}\label{derivee de lenergie}
 \frac{d}{dt}E(u(t))=-a(\|u^{\frac{p}{2}}\nabla u\|_{L^{2}}^{2} - C_{p}\|u\|_{L^{\frac{4}{d}+2+p}}^{\frac{4}{d}+2+p})
 \end{equation}
 and
 \begin{equation}\label{moment}
 \frac{d}{dt}P(u(t))= -2a Im \int\overline{u} |u|^p\nabla u\, .
 \end{equation}
where $C_{p} = \frac{4+2d +pd}{4+2d}$.
\end{remark}
\vspace{0.5cm}

Now we are ready to prove part (4) of Theorem \ref{theoremprincipal}:
\\
Suppose that there exist an initial data $u_0$ with $\|u_0\|_{L^2} \leq \|Q\|_{L^2}$ , such that the corresponding solution $u(t)$ blows up
 with the following law $$\frac{1}{(T-t)^{\beta - \epsilon }} \lesssim \|\nabla u(t)\|_{L^2(\mathbb{R}^{d})} \lesssim \frac{1}{(T-t)^{\beta + \epsilon}},$$
$\text{where}~ \beta;~0 < \beta < \beta(p,d)= \frac{2}{pd}~~ \text{and} ~~0 < \epsilon < \frac{2-\beta pd}{8+pd}$.\\
Recall that
\begin{equation}\label{energie}
\displaystyle{E(u(t)) = E(u_{0}) - a\int_{0}^{t}K(u(\tau))d\tau, \quad t \in [0,T[,}
\end{equation}
where $K(u(t)) =(\|u^{\frac{p}{2}}\nabla u\|_{L^{2}}^{2} - C_{p}\|u\|_{L^{\frac{4}{d}+2+p}}^{\frac{4}{d}+2+p})$.\\
 By  Gagliardo-Nirenberg inequality and (\ref{mass a}), we have:
$$
E(u(t))  \lesssim E(u_0) + \int_0^t\left\| u\right\|_{L^2}^{\frac{4}{d}+p-\frac{p}{2}}\left\|\nabla u\right\|_{L^2}^{2+\frac{pd}{2}}
\lesssim E(u_0) + \int_0^t\left\|\nabla u\right\|_{L^2}^{2+\frac{pd}{2}}
$$
Since the choise of $\epsilon$, we obtain that 
\begin{equation}\label{ksurnabla}
\displaystyle{0 \leqslant \lim_{t \rightarrow T}\frac{\displaystyle{\int_{0}^{t}}
\left\|\nabla u(\tau)\right\|_{L^2}^{2+\frac{pd}{2}}d\tau}{\left\|\nabla u(t)\right\|_{L^2(\mathbb{R}^{d})}^{2}} \lesssim  \lim_{t \rightarrow T} (T-t)^{\frac{1}{2}(2-\beta pd-\epsilon(8+pd))}}=0,
\end{equation}
let
 $$\rho(t) = \frac{\left\|\nabla Q\right\|_{L^2(\mathbb{R}^{d})}}{\left\|\nabla u(t)\right\|_{L^2(\mathbb{R}^{d})}} \quad \text{and} \quad v(t,x)=\rho^{\frac{d}{2}}u(t,\rho x)$$
$(t_{k})_{k}$ a sequence such that $t_{k} \rightarrow T$
and $\rho_{k} = \rho(t_{k}), v_{k} = v(t_{k},.)$. The family $(v_{k})_{k}$ satisfies 

$$\left\|v_{k}\right\|_{L^2(\mathbb{R}^{d})} \leq \left\|u_{0}\right\|_{L^2(\mathbb{R}^{d})}\leq \left\|Q\right\|_{L^2(\mathbb{R}^{d})}
\quad \text{and} \quad \left\|\nabla v_{k}\right\|_{L^2(\mathbb{R}^{d})} = \left\|\nabla Q\right\|_{L^2(\mathbb{R}^{d})}.$$
\bigskip
 Remark that $\displaystyle{\lim_{k \longrightarrow +\infty}} E(v_{k}) =0$, because:
\begin{align}\label{Edevk}
 0\leq \frac{1}{2}(\int |\nabla v_{k}|^2)\bigg(1 - \big(\frac{\int |v_k|^2}{\int |Q|^2}\big)^{2}\bigg)\leq 
E(v_{k}) &= \rho^2_{k}E(u_{0}) -a\rho^2_{k}\int_{0}^{t_{k}}K(u(\tau))d\tau\nonumber\\
&\leq \rho^2_{k}E(u_{0}) + 
\frac{1}{\left\|\nabla u(t_k)\right\|_{L^2(\mathbb{R}^{d})}^{2}}\int_{0}^{t_k}\left\|\nabla u(\tau)\right\|_{L^2}^{2+\frac{pd}{2}}d\tau\nonumber
\end{align}
then using (\ref{ksurnabla}), the energy of $v_k$ tends to $0$.
Which yields 
\begin{equation}\label{vk ver Q}
\displaystyle{\left\|v_{k}\right\|_{L^{\frac{4}{d}+2}}^{\frac{4}{d}+2 } \rightarrow \frac{d+2}{d}\left\|\nabla Q\right\|_{L^2(\mathbb{R}^{d})}^{2}.}
\end{equation}
The family $(v_{k})_{k}$ satisfies the hypotheses of Theorem \ref{limsup} with \\
$$m^{\frac{4}{d}+2} = \frac{d+2}{d}\left\|\nabla Q\right\|_{L^2(\mathbb{R}^{d})}^{2} \quad \text{and} \quad M = 
\left\|\nabla Q\right\|_{L^2(\mathbb{R}^{d})},$$
\bigskip
thus there exists a family $(x_{k})_{k} \subset \mathbb{R}$ and a profile $V \in H^{1}(\mathbb{R})$ with $\left\|V\right\|_{L^2(\mathbb{R}^{d})} \geq \left\|Q\right\|_{L^2(\mathbb{R}^{d})}$, such that,
\begin{equation}\label{convergencefaible}
\displaystyle{\rho^{\frac{d}{2}}_{k}u(t_{k}, \rho_{k}\cdot +  x_{k}) \rightharpoonup V \in H^{1} \quad \text{weakly}.}
\end{equation}
Using (\ref{convergencefaible}), $\forall A \geq 0$
\begin{equation}
 \displaystyle{\liminf_{n\to +\infty}\int_{B(0,A)}\rho_{n}^d|u(t_{n},\rho_{n}x+x_{n})|^{2}dx\geq \int_{B(0,A)}|V|^{2}dx.}\nonumber
 \end{equation}
But $\lim_{n\to +\infty}\frac{1}{\rho_{n}}=+\infty$\,\,thus $\frac{1}{\rho_{n}}> A$, $\rho_{n}A < 1$. This gives immediately:
  
  $$
  \displaystyle{\liminf_{n\to +\infty}\sup_{y\in\mathbb{R}}\int_{|x-y|\leq 1}|u(t_{n},x)|^{2}dx \geq \liminf_{n\to +\infty}
 \int_{|x-x_n| \leq \rho_{n}A}|u(t_{n},x)|^{2}dx \geq \int_{|x|\leq A}|V|^{2}dx.}
  $$

  This it is true for all $A > 0$ thus :
  
  \begin{equation}
  \displaystyle{\liminf_{n\to +\infty}\sup_{y\in\mathbb{R}}\int_{|x-y|\leq 1}|u(t_ {n},x)|^{2}dx\geq \int Q^{2},}
  \end{equation}
  
 then
$$
 \|u_0\|_{L^2} >  \liminf_{n\to +\infty}\sup_{y\in\mathbb{R}}\int_{|x-y|\leq 1}|u(t_{n},x)|^{2}dx \geqslant \|Q\|_{L^2}.
$$
This gives the proof, the fact that $\|u_0\|_{L^2} \leq \|Q\|_{L^2}$.

\section{Global existence.}
In this section, we prove assertion (1) and (2)  of Theorem \ref{theoremprincipal} and Theorem \ref{theorem2}. To prove part (1), we will prove that 
the $H^1$-norm of $u$ is bounded for any time. To prove part (2), we  use generalised 
Gagliardo-Nirenberg inequalities to show that the energy is non increasing.
Finally to prove Theorem \ref{theorem2}, we etablish an a priori estimate on the critical Strichartz norm.

\begin{theorem}
 Let $ p \ge 1  $ for $ d=1,2 $ or $ 1\le p \leq \frac{4}{d-2}$ for $d \ge 3 $, then the initial-value probem (\ref{NLSap}) is locally well posed in $H^{1}(\R^d)$(If $p< \frac{4}{d-2}$ the minimal time
 of the existence depends on $\|u_0\|_{H^1}$.)
.
\end{theorem}
\textbf{Proof:} See \cite{Cazenave1} page 93 Theorem 4.4.1.\\

To prove the following proposition, we will proceed in  the same way as in the section 3.1 in \cite{passota}.
\begin{proposition}\label{global1}
Let $u$ be a solution of (\ref{NLSap}) and $\frac{4}{d-2} \geqslant p > \frac{4}{d}$ then
$$
\|\nabla u(t)\|_{L^2(\mathbb{R}^{d})} \leq \|\nabla u(0)\|_{L^2(\mathbb{R}^{d})}e^{a^{(\frac{-4t}{pd-4})}}.
$$
 
\end{proposition}
\textbf{Proof}:
 Multiply Eq. (\ref{NLSap}) by − $\Delta \overline{u}$ , integrate  and take the imaginary part, this gives
\begin{equation}\label{apriori}
\frac{1}{2}\frac{d}{dt}\int|\nabla u|^2dx + a \int|u|^p|\nabla u |^2 + a\Re\int u\nabla|u|^p\nabla \overline{u}dx= 
-\frac{4}{d}\Im\int u\nabla\overline{u}\Re({u\nabla u})|u|^{\frac{4}{d}-2}.
  \end{equation}

In the l.h.s, a simple calculation shows that the third term rewrites in the form $\frac{p}{4}\int|u|^{p-2}(\nabla|u|^2)^2$.Equation (\ref{apriori})
 becomes:
\begin{equation}\label{apriori1}
 \frac{1}{2}\frac{d}{dt}\int|\nabla u|^2dx + a \int|u|^p|\nabla u |^2 + a\frac{p}{4}\int|u|^{p-2}(\nabla|u|^2)^2  
\leq \frac{2}{d}\int|u|^{\frac{4}{d}}|\nabla u|^2.
\end{equation}
To estimate the r.h.s of (\ref{apriori1}), we rewrite it as ( $p > \frac{4}{d}$)
$$
\int|u|^{\frac{4}{d}}|\nabla u|^2 = \int |u|^{\frac{4}{d}}|\nabla u|^{\frac{8}{pd}}|\nabla u|^{2-\frac{8}{pd}}.
$$
Now by H\"{o}lder inequality we obtain that
$$
\int|u|^{\frac{4}{d}}|\nabla u|^2 \leq (\int|u|^{p}|\nabla u|^{2})^{\frac{4}{pd}}(\int |\nabla u|^{2})^{1-\frac{4}{pd}}.
$$
Then inequality (\ref{apriori1}) takes the form:
$$
\frac{d}{dt}w(t) + 2av(t) \leq \frac{4}{d}v(t)^{\frac{4}{pd}}w(t)^{1-\frac{4}{pd}}.
$$
where $w(t) = \displaystyle{\int} |\nabla u|^2$ and $v(t) = \displaystyle{\int} |u|^{p}|\nabla u|^{2}$.\\
Using Young's inequality $ab \leq \epsilon a^{q} + C\epsilon^{-\frac{1}{q-1}}b^{{q^{\prime}}}$, $\frac{1}{q}+\frac{1}{q^{\prime}}=1$, with
 $q=\frac{pd}{4}$ and $\epsilon = \frac{ad}{2}$ we obtain :

$$\frac{d}{dt}w(t) \leq a^{-\frac{1}{\frac{pd}{4}-1}}w(t).$$
This ensures that:
$$w(t) \leq w(0)e^{a^{\big(-\frac{4t}{pd-4}\big)}}.\hfill{\Box}
$$
This show that the $H^1$-norm of $u$ is bounded for any time and gives directly the proof of part one of Theorem \ref{theoremprincipal} in the case 
 $p>4/d $.\vspace*{5mm}\\
Now we will prove the global existence for  small data, for this we will use the following generalized Gagliardo-Niremberg inequalities (see for instance \cite{Friedman}):
\begin{lemma}\label{Friedman}
 Let $q$, $r$ be  any real numbers satisying $1\leq q$, $ r \leq \infty$, and let $j$, $m$ be any integers satisfying $0\leq j <m$. If $u$ is any functions in 
$C^{m}_{0}(\R^d)$, then 
$$
\|D^j u\|_{L^s} \leq C \|D^m u\|^a_{r} \|u\|_q^{1-a}
$$
where\\ 
$$
\frac{1}{s} = \frac{j}{d} + a(\frac{1}{r}-\frac{m}{d}) + (1-a) \frac{1}{q},
$$
for all $a$ in the interval
$$
\frac{j}{m} \leq a \leq 1,
$$
where $C$ is a constant depending only on $ d$, $m$, $j$,$q$,$r$ and $a$.
\end{lemma}
As a direct consequence we get :
\begin{lemma}\label{termenonlinear}
 Let $1\le p \leq 2$ and $ v \in C^{\infty}_{0}(R^d)$ then:
$$
\int |v|^{\frac{4}{d} + 2 + p} \leq C (\int |\nabla(|v|^{\frac{p+2}{2}})|^2 )\times (\int |v|^{2})^{\frac{2}{d}}.
$$
where $c > 0 $ depending only on $d$ and $p$.
\end{lemma}
\noindent
\textbf{Proof:}
Take $s = \frac{\frac{4}{d} + 2 + p}{1+\frac{p}{2}}$, $q = \frac{2}{1+\frac{p}{2}}$ $r=2$, $j=0$ and $m = 1$,then by Lemma \ref{Friedman} we obtain that:
$$
|u|_{L^{\frac{\frac{4}{d} + 2 + p}{1+\frac{p}{2}}}} \leq C|
\nabla u|^{{\frac{4+2p}{\frac{8}{d} + 4 + 2p}}}_{L^2}|u|^{\frac{\frac{8}{d}}{\frac{8}{d} + 4 + 2p}}_{L^{\frac{2}{1+\frac{p}{2}}}}.
$$
Taking  $u = |v|^{1+\frac{p}{2}}$, we obtain our lemma. \qed\\
Now we can prove the following proposition:
\begin{proposition}\label{ez}
 Let $1\le p \leq 2$. There exists  $ 0 < \alpha = \alpha(p,d)< \|Q\|_{L^2}$ , such that  for any $u_0 \in H^1$ with $\|u_0\|_{L^2} < \alpha$, it holds
$$
\frac{d}{dt}E(u(t)) \leq 0, ~~\forall t > 0.
$$
\end{proposition}
\noindent
\textbf{Proof:}
We can whrite that:
$$
\frac{d}{dt}E(u(t)) = a \Bigl( C_p\int  |u |^{\frac{4}{d}+p+2} - \frac{4}{(p+2)^2}\int  |\nabla (|u|^{\frac{p+2}{2}}) |^2\Bigr),
$$

then by Lemma \ref{termenonlinear} we obtain that:
$$
\frac{d}{dt}E(u(t)) \leq a(\int  |\nabla (|u|^{\frac{p+2}{2}}) |^2)(C_pC(\int |u|^{2})^{\frac{2}{d}} -\frac{4}{(p+2)^2} )
$$
 Choosing $\alpha^{\frac{2}{d}} < \frac{4}{(p+2)^2}\frac{1}{C_pC}$, and using  that $ \|u(t)\|_{L^2} \le \|u_0\|_{L^2} $ for all $ t\ge 0 $ (see \ref{mass a} below)  we get the result. \qed\\
Now the proof of  part (2) of Theorem \ref{theoremprincipal} follows from  the sharp Gagliardo-Nirenberg inequality :
\begin{equation}
\forall u \in H^1,  E(u) \geq \frac{1}{2}(\int |\nabla u|^2)\bigg(1 - \big(\frac{\int |u|^2}{\int |Q|^2}\big)^{\frac{2}{d}}\bigg). \nonumber
\end{equation}
Proposition \ref{ez} together with the above inequality ensure that the $ H^1$-norm of $ u $ is uniformly bounded in time. This leads to the   global existence result  for small initial data   when  $ 1\le p\le 2  $.  

\subsection{Critical case ( $p = \frac{4}{d}$)}
Now we will treat the critical case and prove Theorem \ref{theorem2}. First let us  prove that, if the solution blows up in finite time $T$, then
 \textbf{$\|u\|_{L^{\frac{4}{d}+2}([0,T];L^{\frac{4}{d}+2}(R^d))} = +\infty$}.
\begin{proposition}\label{propessentiel}
 Let $u$ be the unique maximal solution of (\ref{NLSap}) in $[0,T^{*})$; if $T^{*} < \infty$, then $\|u\|_{L^{\sigma}([0,T],L^{\sigma})} = \infty$ 
where $\sigma = \frac{4}{d} + 2$.
\end{proposition}
To prove this claim, denoting by $ S(\cdot) $  the free evolution of the linear Schr\"odinger equation and defining   the notion of admissible pair  in the following way :  An ordered pair $(q, r )$ is called admissible if $\frac{2}{q} + \frac{d}{r} = \frac{d}{2}$, $2 < q \leq \infty $, we will use the following proposition:
\begin{proposition}\label{propodeCazenave}
 There exists $\delta > 0$ with the following property. If $u_0 \in  L^2(R^ d)$ and $ T \in (0,\infty]$ are such
that $\|S(.)u_0\|_{L^{\sigma}([0,T],L^{\sigma})} < \delta$, there exists a unique solution $u \in C([0,T],L^2(R^d))\cap L^{\sigma}([0,T], L^{\sigma}(R^n))$ of 
(\ref{NLSap}).
In addition, $u\in L^q([0,T],L^r(R^d))$ for every admissible pair $(q,r)$; for $t \in [0,T]$.Finally, $u$
depends continuously in $C([0,T],L^2(R^d)) \cap L^{\sigma}([0,T],L^{\sigma}(R^d))$ on $u_0 \in L^2(R^n)$. If $u_0 \in  H^1(R^d)$, then
$u \in C([0,T],H^1(R^d))$.
\end{proposition}
See \cite {Cazenave2} for the proof.\\ 
\\
We need also the following lemma ( see \cite{Cazenave2}):
\begin{lemma}\label{lemmaCazenave}
 Let $T \in (0,\infty]$, let $\sigma=\frac{4}{d}+2$, and let $(q,r)$ be an admissible pair. Then,
whenever $u \in L^{\sigma}([0,T],L^{\sigma}(\R^d))$, it follows that $F(u)=-i\displaystyle{\int_0^t S(t-s)(|u|^{\frac{4}{d}}u + ia|u|^{p}u)ds }\in C([0,T],H^{-1}(\R^d))\cap L^q(0,T,L^r(\R^d))$. Furthermore,
there exists $K$, independent of $T$, such that
\begin{equation}\label{F}
\|Fv - Fu\|_{L^q(]0,T[,L^r)} < K ( \|u\|^{\frac{4}{d}}_{L^{\sigma}(]0,T[,L^{\sigma}(\R^d))}+ \|v\|^{\frac{4}{d}}_{L^{\sigma}(]0,T[,L^{\sigma}(\R^d))} )
 \|u-v\|_{L^{\sigma}(]0,T[,L^{\sigma}(\R^d))}
\end{equation}
for every $u$, $v\in L^{\sigma}(]0,T[,L^{\sigma}(\R^d)).$
\end{lemma}
\noindent
\textbf{Proof of Proposition \ref{propessentiel}:}

Let $u_0 \in L^2(\R^d)$. Observe that$ \|S(.)u_0\|_{L^{\sigma}(0,T.L^{\sigma})} \longrightarrow 0$ as $T \longrightarrow  0$. Thus for
sufficiently small $T$, the hypotheses of Proposition \ref{propodeCazenave} are satisfied. Applying iteratively this proposition,
we can construct the maximal solution $u \in C([0,T^*),L^2(\R^d)))\cap L^{\sigma} ([0,T^*),L^{\sigma}(\R^d))$ of (\ref{NLSap}). We proceed by contradiction,
 assuming that  
$T^*< \infty$, and $\|u\|_{L^{\sigma}(]0,T[,L^{\sigma})} < \infty$. Let $t \in[0,T^*)$. For every $s \in [0,T^{*}-t)$ we have
$$
S(s)u(t) = u(t+s) - F(u(t+\cdotp))(s).
$$
From  (\ref{F}), we thus obtain
$$
\|S(.)u(t)\|_{ L^{\sigma} ([0,T^*-t),L^{\sigma}(\R^d))} \leq \|u\|_{L^{\sigma}(]t,T^{*}[,L^{\sigma})} + K(\|u\|_{L^{\sigma}(]t,T^*[,L^{\sigma})})^{\frac{4}{d}+1}
$$
Therefore, for $t$ fixed close enough to $T^*$, it follows that
$$
\|S(.)u(t)\|_{ L^{\sigma} ([0,T^*-t),L^{\sigma}(\R^d))} \leq \delta.
$$
Applying Proposition \ref{propodeCazenave}, we find that $u$ can be extended after $T^*$, which contradicts the maximality.

\begin{corollary}
 For $p = \frac{4}{d}$, the solution of (\ref{NLSap}) is global.
\end{corollary}
\textbf{Proof}
Multiply equation (\ref{NLSap}) by $\overline{u}$, and take the imaginary part to obtain:
$$
\frac{d}{dt}\|u(t)\|^2_{L^2} + 2a\|u\|^{\frac{4}{d+2}}_{L^{\frac{4}{d}+2}}=0.
$$ 
Hence $\forall t \in \R_{+}$
$$
 \|u\|_{L^{\frac{4}{d}+2}[0,t[L^{\frac{4}{d}+2}(R^d)} \leq \frac{1}{2a} \|u_0\|^2_{L^2}.
$$
\\
 The global existence follows then directly from Propostion \ref{propessentiel}. Now to finish the proof of Theorem \ref{theorem2}, 
we will prove the scattering:\\
Let $v(t) = e^{-it\Delta}u(t) :=S(-t)u(t)$ then
\begin{equation}
 v(t) = u_{0} + i\int_{0}^{t}S(-s)(|u(s)|^{\frac{4}{d}}u(s) + ia|u(s)|^{\frac{4}{d}}u)ds. \nonumber
\end{equation} 
Therefore for $0 < t <\tau$,
$$ 
v(t)-v(\tau) = i\int_{\tau}^{t}S(-s)(|u(s)|^{\frac{4}{d}}u(s) + ia|u(s)|^{\frac{4}{d}}u)ds.
$$
It follows from Strichartz's estimates (see the proof of Lemma \ref{us1}) that:
$$
\|v(t)-v(\tau)\|_{L^2} = \|i\int_{\tau}^{t}S(-s)(|u(s)|^{\frac{4}{d}}u(s) + ia|u(s)|^{\frac{4}{d}}u)ds\|_{L^2} 
\leq C \parallel u\parallel_{L^{\frac{4}{d}+2}([t,\tau]\times\mathbb{R}^{d})}^{\frac{4}{d}+1}.
$$
The right hand side goes to zero when $t, \tau \longrightarrow + \infty$,  then scattering follows from the Cauchy criterion.\\
 This completes the proof of Theorem \ref{theorem2}.

\section{Blow up solution.}
In this section, we will prove the existence of the explosive solutions in the case $ 1\le p<4/d$. 
\begin{theorem}\label{theorem 3} Let $ 1\le p <4/d $. 
 There exist  a set of initial data $\Omega$ in $H^1$, such that  for any $0 < a < a_0$ with $a_0=a_0(p)$ small enough, the emanating solution $u(t)$  to (\ref{NLSap}) blows up in finite time in the log-log regime.
\end{theorem}
The set of initial data $\Omega$  is the set described  in \cite{MerleRaphael1} �in order to initialize the log-log regime. It is open in $H^1$. Using the continuity  with regard to the initial data and the parameters, we easily obtain the following corollary:
\begin{corollary}\label{colloraireperturbation}
 Let $ 1\le p <4/d $ and $u_{0} \in H^1$ be an initial data  such that the corresponding solution $u(t)$ of (\ref{NLS}) blows up in the loglog regime. There exist $\beta_{0} > 0$ and $a_{0} > 0 $ such that if $v_{0} = u_{0} + h_{0}$, $\left\|h_{0}\right\|_{H^{1}} \leq \beta_{0}$ and $a \leq a_{0}$, the solution $v(t)$ for (\ref{NLSap}) with the initial data $v_{0}$ blows up in finite time.
\end{corollary}
Assertion (3) of Theorem \ref{theoremprincipal} now follows directly from  this corollary together the results of  \cite{MerleRaphael1} on the $ L^2$-critical NLS equation and a  scaling argument in order to drop the smallness  condition on the damped coefficient $ a>0$. \\
\noindent 
Now to prove Theorem \ref{theorem 3}, 
we look for a solution of (\ref{NLSap}) such that for $t$ close enough to blowup time, we shall have the following decomposition:
\begin{equation}\label{decomposition}
u(t,x)=\frac{1}{\lambda^{\frac{d}{2}}(t)}(Q_{b(t)} + \epsilon)(t,\frac{x-x(t)}{\lambda(t)})e^{i\gamma(t)},
\end{equation}

for some geometrical parameters $(b(t),\lambda(t), x(t),\gamma(t)) \in (0,\infty)\times(0,\infty)\times\mathbb{R}^{d}\times\mathbb{R}$, here $\lambda(t)\sim \frac{1}{\|\nabla u(t)\|_{L^2}}$,
 and the profiles $Q_{b}$ are suitable deformations of $Q$ related to some extra degeneracy
of the problem.\\

 Note that  we will abbreviated our proof because it is very very close to the case of  
linear damping ($p=0$ see Darwich\cite{Darwich}). Actually, as noticed in \cite{Planchon1}, we only need to prove that in the log-log regime
the $L^2$ norm does not grow, and the growth of the energy( resp the momentum) is below $\frac{1}{\lambda^2}$ (resp $\frac{1}{\lambda}$) .
In this paper, we will prove that in the
log-log regime, the growth of the energy and the momentum are bounded by:
$$E(u(t))\lesssim log (\lambda(t))\lambda(t)^{-\frac{pd}{2}}, ~~ 	P(u(t)) \lesssim log(\lambda(t))\lambda(t)^{1-\frac{pd}{4}}.$$
Let us recall that a fonction u :$[0,T]\longmapsto H^1$ follows the log-log regime if the following
 uniform controls on the decomposition (\ref{decomposition}) hold on $[0,T]$:
\begin{itemize}
\item {Control of $b(t)$}
\begin{equation}\label{b petit}
 b(t) > 0 , ~ b(t) < 10 b(0).
\end{equation}

\item Control of $\lambda$:
\begin{equation}\label{control of lambda}
\lambda(t) \leq e^{-{e^\frac{\pi}{100b(t)}}}
\end{equation}
and the monotonicity of $\lambda$:
\begin{equation}\label{monotonicity}
\lambda(t_2) \leq\frac{3}{2} \lambda(t_{1}), \forall~0\leq t_{1} \leq t_{2} \leq T.
\end{equation}
Let $k_{0} \leq k_{+} $ be integers and $T^{+} \in [0,T]$  such that
\begin{equation}\label{lambda 0 et lambda T}
\frac{1}{2^{k_{0}}} \leq \lambda(0) \leq \frac{1}{2^{k_{0}-1}}, \frac{1}{2^{k_{+}}} \leq \lambda(T^{+}) \leq \frac{1}{2^{k_{+}-1}}
\end{equation}
and for $k_{0} \leq k \leq k_{+} $, let $t_{k}$ be a time such that
\begin{equation}
\lambda(t_{k}) = \frac{1}{2^k},
\end{equation}
then we assume the control of the doubling time interval:
\begin{equation}\label{tk}
t_{k+1} - t_{k} \leq k \lambda^2(t_{k}).
\end{equation}
\item control of the excess of mass:
\begin{equation}\label{controlepsilon}
\int\left|\nabla \epsilon(t)\right|^2 + \int\left|\epsilon(t)\right|^2e^{-\left|y\right|} \leq \Gamma_{b(t)}^{\frac{1}{4}}.
\end{equation}
\end{itemize}
\bigskip

\subsection{Control of the energy and the kinetic momentum in the log-log regime}
\bigskip
We recall the Strichartz estimates. An ordered pair $(q, r )$ is called admissible if $\frac{2}{q} + \frac{d}{r} = \frac{d}{2}$, $2 < q \leq \infty $.
 We define the Strichartz norm of functions $u : [0, T]\times \R^{d} \longmapsto C$ by:
\begin{equation}\label{us0}
\left\|u\right\|_{S^{0}([0,T]\times\mathbb{R}^{d})} = \sup_{(q,r) admissible}\left\|u\right\|_{L^{q}_{t}L^{r}_{x}([0,T]\times\mathbb{R}^{d})}
\end{equation}
and
\begin{equation}\label{us1}
\left\|u\right\|_{S^{1}([0,T]\times\mathbb{R}^{d})} =  \sup_{(q,r) admissible}\left\|\nabla u\right\|_{L^{q}_{t}L^{r}_{x}([0,T]\times\mathbb{R}^{d})}
\end{equation}
We will sometimes abbreviate $S^i([0,T]\times \mathbb{R}^{2})$ with $S^{i}_{T}$  or $S^{i}[0, T]$, $i= 1,2$. Let us denote the H\"older dual exponent of $q$ by $q^{\prime}$ so that $\frac{1}{q} + \frac{1}{q^{\prime}} = 1$. The Strichartz estimates may be
expressed as:
\begin{equation}\label{estimationschwartz}
\left\|u\right\|_{S^{0}_{T}} \leq \left\|u_0\right\|_{L^{2}} + \left\|(i\partial_{t} + \Delta )u\right\|_{L^{q^{\prime}}_{t}L^{r^{\prime}}_{x}}
\end{equation}
where $(q, r )$ is any admissible pair.
Now we will derive an estimate on the energy, to check that it remains small with respect to $\lambda^{-2}$:

\begin{lemma}\label{control energie}
Assuming that (\ref{b petit})-(\ref{controlepsilon}) hold, then
 the energy and kinetic momentum of the solution $u$ to (\ref{NLSap}) are controlled on $[0, T]$ by:\\
\begin{equation}\label{control de lenergie}
|E(u(t))|\leq C(log(\lambda(t))\lambda(t)^{-\frac{pd}{4}}),
\end{equation}

\begin{equation}\label{control de moment 1}
|P(u(t))|\leq C (log(\lambda(t))\lambda(t)^{1-\frac{pd}{4}}).
\end{equation}
\end{lemma}
To prove this lemma, we shall need the following one:\\
\begin{lemma}\label{us1} 
Let $u$ be a solution of $(\ref{NLSap})$ emanating for $u_{0}$ in $H^1$. Then $u$ $\in$ $C([0,\Delta T], H^1)$ where $\Delta T = \left\|u_{0}\right\|^{\frac{d-4}{d}}_{L^2}\left\|u_{0}\right\|_{H^{1}}^{-2}$, and we have the following control\\
\centerline{
$\left\|u\right\|_{{S^{0}[t,t+\Delta T]}} \leq 2\left\|u_{0}\right\|_{L^2}$ , $\left\|u\right\|_{{S^{1}[t,t+\Delta T]}} \leq 2 \left\|u_{0}\right\|_{H^1(\mathbb{R}^{d})}.$}
\end{lemma}
Proof of \textbf{Lemma \ref{control energie}:}
According to (\ref{tk}) each interval $[t_k, t_{k+1}]$, can be divided into $k$ intervals, $[\tau_{k}^{j},\tau_{k}^{j+1}]$  
such that the estimates of the previous lemma are true.
From (\ref{derivee de lenergie}) and the  Gagliardo-Nirenberg inequality, we obtain that:
$$
\frac{d}{dt}E(u(t)) \lesssim \left\| u\right\|_{L^2}^{\frac{4}{d}+p-\frac{pd}{2}}\left\|\nabla u\right\|_{L^2}^{2+\frac{pd}{2}}
$$
Using (\ref{mass a}) this gives
$$
\int_{\tau_{k}^{j}}^{\tau_{{k}^{j+1}}}\frac{d}{dt}E(u(t))dt \leq C\int_{\tau_{k^j}}^{\tau_k^{j+1}}\left\|\nabla u(t)\right\|_{L^2}^{2+\frac{pd}{2}},
$$
then by Lemma \ref{us1}, we obtain that:
$$\int_{\tau_{k}^{j}}^{\tau_{{k}^{j+1}}}\frac{d}{dt}E(u(t))dt  \leq C (\tau_{k^{j+1}}-\tau_{k^j})\lambda^{-2-\frac{pd}{2}}(\tau_{k^j})$$
Note that $\tau_{k}^{j+1}-\tau_{k}^{j} \sim \lambda^2(\tau_{k}^{j}) \sim \lambda^2(t_k)$, then
$$
\int_{\tau_{k}^j}^{\tau_{k}^{j+1}}\frac{d}{dt}E(u(t))dt \leq C \lambda^{-\frac{pd}{2}}(t_{k})
$$
 Summing from $j=1$ to $J_k \leq C K$, we obtain that:
$$
\sum_{j=1}^{J_k}\int_{\tau_{k}^j}^{\tau_{k}^{j+1}}\frac{d}{dt}E(u(t))dt \leq C k\lambda^{-\frac{pd}{2}}(t_{k})
$$
Now taking $T^+ = T$ and summing from $K_0$ to $K^+$, we obtain:
$$
\int_0^{T^+}\frac{d}{dt}E(u(t))dt\leq C K^+\lambda^{-\frac{pd}{2}}(T^+)\lesssim C log (\lambda(T))\lambda^{-\frac{pd}{2}}(T).
$$
Note that log $(\lambda(T))\lambda^{-\frac{pd}{2}}(T)$ is  small with to respect $\frac{1}{\lambda^2}$ because $p < \frac{4}{d}$.\\
\\
Now we prove \textbf{(\ref{control de moment 1})}:
From (\ref{moment}) we have :

$$
|\frac{d}{dt}P(u(t))| \leq \int |\overline{u}||\nabla u||u|^p
$$
By  Gagliardo-Nirenberg inequality we have:
$$
\left\| u\right\|^{2p+2}_{L^{2p+2}} \leq  \left\| u\right\|_{L^2}^{2p+2-dp}\left\| \nabla u\right\|_{L^2}^{dp}
$$
then
$$
\frac{d}{dt}P(u(t)) \leq (\int |u|^{2(p+1)})^{\frac{1}{2}}\left\| \nabla u\right\|_{L^2} \leq C \left\| \nabla u\right\|^{1+\frac{pd}{2}}_{L^2}.
$$
Then:
$$
\int_{\tau_{k}^j}^{\tau_{k}^{j+1}}\frac{d}{dt}P(u(t))\leq C (\tau_{k}^{j+1}-\tau_{k}^{j})\left\| \nabla u(\tau_{k}^j)\right\|^{1+\frac{pd}{2}}_{L^2}
\leq C\left\| \nabla u(t_k)\right\|^{-1+\frac{pd}{2}}_{L^2}
$$
Summing successively into $j$ and $k$ we obtain that:$$
\int_0^{T^+}\frac{d}{dt}P(u(t)) \lesssim \text{log} (\lambda(T^+))\lambda^{1-\frac{pd}{2}}(T^+).$$
Remark that this quantity is small with to respect $\frac{1}{\lambda}$ because $p< \frac{4}{d}$.$\hfill{\Box}$



\begin{thebibliography}{10}
\bibitem{Sparber}
P. Antonelli and C. Sparber.\newblock{ \it Global well-posedness for cubic NLS with nonlinear damping}.
\newblock{ Comm. Partial Differential Equations}, 35 (2010) 4832–4845.


\bibitem{Ber}
H.~Berestycki and P.-L. Lions.
\newblock{\it Nonlinear scalar field equations. {II}. {E}xistence of infinitely
  many solutions}.
\newblock { Arch. Rational Mech. Anal.}, 82(1983):347--375.

\bibitem{Cazenave1}
T.~Cazenave.
\newblock {\it Semilinear {S}chr\"odinger equations}, volume~10 of {\em Courant
  Lecture Notes in Mathematics}.
\newblock New York University Courant Institute of Mathematical Sciences, New
  York, 2003.
\bibitem{Cazenave2}
T.~Cazenave and F.~Weissler.
 \newblock{\it Some remarks on the nonlinear Schrödinger equation in the subcritical case}.
 \newblock New methods and results in nonlinear field equations (Bielefeld, 1987), 59–69, Lecture Notes in Phys., 347, Springer, Berlin, 1989.

\bibitem{Colliander}
J.~Colliander and P.~Raphael.
\newblock {\it Rough blowup solutions to the {$L^2$} critical {NLS}}.
\newblock { Math. Ann.}, 345(2009):307--366.
\bibitem{Darwich}
M.~Darwich.
\newblock{\it Blowup for the Damped $L^2$critical nonlinear Shr\"odinger equations.}
\newblock{Advances in Differential Equations.}
volume 17, Numbers 3-4 (2012),337-367. 
\bibitem{Fibich}
G.~Fibich and F.~Merle.
\newblock {\it Self-focusing on bounded domains}.
\newblock { Phys. D}, 155(2001):132--158.
\bibitem{Fibich1}
G.~ Fibich and M.~ Klein. \newblock{\it Nonlinear-damping continuation of the nonlinear Schrödinger equation-a numerical study}.
\newblock{ Physica D}, 241 (2012), 519-527.

\bibitem{Friedman}
A.~Friedman
\newblock{\it Partial Differential Equations.}
\bibitem{Hmidi}
T.~Hmidi and S.~Keraani.
\newblock{\it Blowup theory for the critical nonlinear {S}chr\"odinger equations
  revisited}.
\newblock {Int. Math. Res. Not.}, 46(2005):2815--2828.

\bibitem{Kato}
T.~Kato.
\newblock {\it On nonlinear {S}chr\"odinger equations}
\newblock {Ann. Inst. H. Poincar\'e Phys. Th\'eor.}, 46(1987):113--129.

\bibitem{Kwong1}
M.K Kwong.
\newblock {\it Uniqueness of positive solutions of {$\Delta u-u+u^p=0$} in {${\bf
  R}^n$}}.
\newblock { Arch. Rational Mech. Anal.}, 105(1989):243--266.

\bibitem{Lions1}
P.-L. Lions.
\newblock {\it The concentration-compactness principle in the calculus of
  variations. {T}he locally compact case. {II}}.
\newblock { Ann. Inst. H. Poincar\'e Anal. Non Lin\'eaire}, 1(1984):223--283.
\bibitem{Merleseul}
F. Merle. Determination of blow-up solutions with minimal mass for
nonlinear Schr \"{o}inger equations with critical power, Duke Math. J. 69:2,
(1993), 427-454.

\bibitem{MerleRaphael1}
F.~Merle and P.~Raphael.
\newblock {\it Blow up dynamic and upper bound on the blow up rate for critical
  nonlinear {S}chr\"odinger equation}.
\newblock In { Journ\'ees ``\'{E}quations aux {D}\'eriv\'ees {P}artielles''
  ({F}orges-les-{E}aux, 2002)}, pages Exp. No. XII, 5. Univ. Nantes, Nantes,
  2002.

\bibitem{Merle2}
F.~Merle and P.~Raphael.
\newblock {\it Sharp upper bound on the blow-up rate for the critical nonlinear
  {S}chr\"odinger equation}.
\newblock { Geom. Funct. Anal.}, 13(2003):591--642.

\bibitem{Merle3}
F.~Merle and P.~Raphael.
\newblock {\it On universality of blow-up profile for {$L^2$} critical nonlinear
  {S}chr\"odinger equation}.
\newblock { Invent. Math.}, 156(2004):565--672.

\bibitem{Merle5}
F.~Merle and P.~Raphael.
\newblock {\it Profiles and quantization of the blow up mass for critical nonlinear
  {S}chr\"odinger equation}.
\newblock {Comm. Math. Phys.}, 253(2005):675--704.

\bibitem{Merle4}
F.~Merle and P.~Raphael.
\newblock {\it On a sharp lower bound on the blow-up rate for the {$L^2$} critical
  nonlinear {S}chr\"odinger equation}.
\newblock { J. Amer. Math. Soc.}, 19(2006):37--90 (electronic).

\bibitem{ohta}
M.~Ohta and G.~Todorova.
\newblock {\it Remarks on global existence and blowup for damped nonlinear
  {S}chr\"odinger equations}.
\newblock { Discrete Contin. Dyn. Syst.}, 23(2009):1313--1325.

\bibitem{passota}
T. Passota, C. Sulemb and P.L. Sulem.
\newblock{\it Linear versus nonlinear dissipation for critical NLS equation.}
\newblock{Physica D}, 203 (2005) 167–184


\bibitem{Planchon1}
F.~Planchon and P.~Rapha{\"e}l.
\newblock {\it Existence and stability of the log-log blow-up dynamics for the
  {$L^2$}-critical nonlinear {S}chr\"odinger equation in a domain}.
\newblock {Ann. Henri Poincar\'e}, 8(2007):1177--1219.

\bibitem{Merle6}
P.~Raphael.
\newblock {\it Stability of the log-log bound for blow up solutions to the critical
  non linear {S}chr\"odinger equation}.
\newblock { Math. Ann.}, 331(2005):577--609.

\bibitem{Tsutsumi}
M.~Tsutsumi.
\newblock {\it Nonexistence of global solutions to the {C}auchy problem for the
  damped nonlinear {S}chr\"odinger equations}.
\newblock { SIAM J. Math. Anal.}, 15(1984):357--366.

\bibitem{weinst}
M.I. Weinstein.
\newblock {\it Nonlinear {S}chr\"odinger equations and sharp interpolation
  estimates}.
\newblock { Comm. Math. Phys.}, 87(1982/83):567--576.

\end{thebibliography}
\end{document}